\documentclass[a4paper,12pt]{article}

\usepackage{amsfonts}
\usepackage{amscd,color}
\usepackage{amsmath,amsfonts,amssymb,amscd}
\usepackage{indentfirst,graphicx,epsfig}
\usepackage{graphicx,psfrag}
\input{epsf}
\usepackage{graphicx}
\usepackage{epstopdf}
\usepackage{caption}

\newtheorem{thm}{Theorem}[section]

\newtheorem{conjecture}[thm]{Conjecture}
\newtheorem{lem}[thm]{Lemma}

\def\qed{\hfill \nopagebreak\rule{5pt}{8pt}}

\title{\bf Solutions for two conjectures on\\[2mm] kaleidoscopic
edge-colorings\footnote{Supported by NSFC No.11371205, 11531011.}}

\author{\small Xueliang Li, Xiaoyu Zhu \\
{\small  Center for Combinatorics and LPMC}\\
{\small Nankai University, Tianjin 300071, P.R. China}\\
{\small Email: lxl@nankai.edu.cn;\ zhuxy@mail.nankai.edu.cn}\\
}
\date{}

\begin{document}
\maketitle
\begin{abstract}
For an $r$-regular graph $G$, we define an edge-coloring $c$ with colors from $\{1,2,\cdots,$ $k\}$, in such a way that any vertex of $G$ is incident to at least one edge of each color. The multiset-color $c_m(v)$ of a vertex $v$ is defined as the ordered tuple $(a_1,a_2,\cdots ,a_k)$, where $a_i \ (1\leq i\leq k)$ denotes the number of edges with color $i$ which are incident with $v$ in $G$. Then this edge-coloring $c$ is called a {\it $k$-kaleidoscopic coloring} of $G$ if every two distinct vertices in $G$ have different multiset-colors and in this way the graph $G$ is defined as a {\it $k$-kaleidoscope}. In this paper, we determine the integer $k$ for a complete graph $K_n$ to be a $k$-kaleidoscope, and hence solve a conjecture in [P. Zhang, A Kaleidoscopic View of Graph Colorings, Springer, New York, 2016] that for any integers $n$ and $k$ with $n\geq k+3 \geq 6$, the complete graph $K_n$ is a $k$-kaleidoscope. Then, we  construct an $r$-regular $3$-kaleidoscope of order $\binom{r-1}{2}-1$ for each integer $r\geq 7$, where $r\equiv 3\ (\text{mod}\ 4)$, which solves another conjecture in the same book on the maximum order for $r$-regular $3$-kaleidoscopes.
\vspace{5 mm}

\vspace{5 mm}

\end{abstract}

\section{Introduction}

In this paper, all graphs are simple, undirected and finite. For notation and terminology we follow the book \cite {BM}. An $edge\text{-}coloring$ for a graph $G$ is a mapping from the edges of $G$ to a finite number of colors. In the early days, many classical colorings were put forward and studied such as proper edge-coloring, list edge-coloring, acyclic edge-coloring and so on. Recently, based on a variety of application instances in different fields, a number of new edge-colorings were put forward. For example, the rainbow edge-coloring has received wide attention due to its close connection with network security and many valuable results were derived in the papers such as \cite{CJMZ, LSS} and the book \cite{LS}.

In this paper, we want to specify more about another kind of edge-coloring, that is, the $kaleidoscopic \ coloring$, which was first introduced in \cite{CEZ}. Assume that a group of $n$ computers, each has $r$ ports on the back, are needed to build a network. There are $k$ kinds of connections altogether and every two distinct computers can build at most one connection between them. It is necessary to use every port so that the fail-safe connections would be maximized. Furthermore, distinct computers must have different numbers of types of connections so that the computer engineer is able to distinguish them. The above fact is an application instance for kaleidoscopic coloring. Actually it can model a lot of situations and can be applied to many fields such as computer science and telecommunications. Next we will give the definition of kaleidoscopic coloring.

For an $r$-regular graph $G$, we define an edge-coloring to its edges with the colors $[k]=\{1,2,3,\cdots,k\} \ (k\geq3)$ such that any vertex in $G$ is incident to at least one edge of each color. For a color set $S=\{i_1,i_2,\cdots,i_s\}$, the $S$-tuple of a vertex $v$ is defined as $(a_{i_1},a_{i_2},\cdots,a_{i_s})$, where $a_{i_j}\ (1\leq j\leq s)$ denotes the number of edges with color $i_j$ which are incident with $v$ in $G$. In particular, the {\it multiset-color $c_m(v)$} of the vertex $v$ is an $S$-tuple for $S=[k]$. Then this edge-coloring $c$ is called a {\it $k$-kaleidoscopic coloring} of $G$ if every two distinct vertices in $G$ have different multiset-colors and in this way the graph $G$ is a $k$-$kaleidoscope$.

As one knows that a proper edge-coloring of a graph $G$ is to factorize the edge set of $G$ into   $\mathcal{F}=\{F_1,F_2,\cdots, F_k\}$ such that $F_i\ (1\leq i\leq k)$ is an independent edge set, what we need for the kaleidoscopic coloring is a factorization that enables distinct vertices have distinct $[k]$-tuples, which is defined as $(\text{deg}_{F_1}(v),\text{deg}_{F_2}(v),\cdots,\text{deg}_{F_k}(v))$ for the vertex $v$ in $G$.

It is well known that every connected graph has at least two vertices with the same degree. Actually there are exactly one connected graph $G$ of order $n$ containing only two vertices of the same degree. We then describe the graph as follows. Label the vertices of $G$ as $v_1,v_2,\cdots ,v_n$, and add an edge $v_iv_j$ if and only if $i+j\geq n+1$, it is obvious that $v_{\lfloor{\frac{n}{2}}\rfloor}$ and $v_{\lfloor{\frac{n}{2}}\rfloor+1}$ share the same degree of $\lfloor{\frac{n}{2}}\rfloor$. Therefore for any decomposition of any $r$-regular graph $G$ into two graphs $D_1$ and $D_2$, at least two vertices $u$ and $v$ have the property: deg$_{D_1}$(u)=deg$_{D_1}$(v) and deg$_{D_2}$(u)=deg$_{D_2}$(v). So any $r$-regular graph $G$ is definitely not a $2$-kaleidoscope, thus our discussion starts with the integer $3$. It can be easily seen that $r>k$ is required for an $r$-regular graph to be a $k$-kaleidoscope according to the definition. But in fact $r>k+1$
holds since $r=k+1$ would imply that at most $k$ different $[k]$-tuples satisfy the demand of at least $r+1$ distinct $[k]$-tuples for the vertices in $G$, a contradiction. Combining with the above, we come to the conclusion that for an $r$-regular graph $G$ of order $n$, the integer $k$ for $G$ to be a $k$-kaleidoscope can only be integers chosen from $3$ to $n-3$.

In \cite{Z}, the author solved the cases for $k$ to be $3$ or $n-3$ when $G$ is complete. And they also posed the following conjecture.

\begin{conjecture}\cite{Z}\label{conjecture1.1}
For integers $n$ and $k$ with $n\geq k+3\geq 6$, the complete graph $K_n$ is a $k$-kaleidoscope.
\end{conjecture}

Another concerning problem in \cite{Z} is the maximum order for the $r$-regular $3$-kaleidos-
copes. A simple calculation shows that there are $\binom{r-1}{2}$ different $[3]$-tuples altogether for an $r$-regular graph. So the number of vertices in an $r$-regular graph can not exceed $\binom{r-1}{2}$ as long as it is a $3$-kaleidoscope according to the definition. Since for the integer $r\equiv 3\ (\text{mod} \ 4)$, $\binom{r-1}{2}$ is odd, so there is no $r$-regular graph with order $\binom{r-1}{2}$. Thus the largest possible order for an $r$-regular $3$-kaleidoscope is $\binom{r-1}{2}-1$. The author in \cite{Z} proved that for any $r\ (r\geq5)$ such that $r\not\equiv3\ (\text{mod}\ 4)$, there exists an $r$-regular $3$-kaleidoscope of order $\binom{r-1}{2}$. Furthermore, the following conjecture was posed in the same book.

\begin{conjecture}\cite{Z}\label{conjecture1.2}
For each integer $r\geq7$ where $r\equiv3\ (\text{mod} \ 4)$, there is an $r$-regular $3$-kaleidoscope of order $\binom{r-1}{2}-1$.
\end{conjecture}

In this paper, we solve these two conjectures, and give their proofs respectively in Sections $2$ and $3$.

\section{Proof of Conjecture \ref{conjecture1.1}}

Before we give the proof of Conjecture \ref{conjecture1.1}, two auxiliary lemmas are stated as follows.

\begin{lem}\cite{Z}\label{lemma2.1}
For each integer $n\geq6$, the complete graph $K_n$ is a $3$-kaleidoscope.
\end{lem}
\begin{lem}\cite{Z}\label{lemma2.2}
For each integer $k\geq3$, the complete graph $K_{k+3}$ is a $k$-kaleidoscope.
\end{lem}

\textbf{Proof of Conjecture \ref{conjecture1.1}:} The case when $k=3$ is verified in Lemma \ref{lemma2.1}. For the complete graph $K_n$, we give our proof by induction on $n$. The case when $n=6$ can easily be verified that it satisfies Conjecture \ref{conjecture1.1} according to Lemma \ref{lemma2.1}. We suppose that $K_m\ (m\geq6)$ is a $k$-kaleidoscope for any $k\ (3\leq k\leq m-3)$, where $m$ is any integer smaller than $n$. We distinguish two cases to clarify.

\textbf{Case 1.} When $\lceil\frac{n}{2}\rceil\leq k\leq n-3$, we all know that for $K_m \ (m\geq4 \ is\ even)$, $K_m$ can be decomposed into $\frac{m}{2}-1$ Hamiltonian cycles $H_1,H_2,\cdots,H_{\frac{m}{2}-1}$ and a perfect matching $F$. We then put the $(n-3)$-kaleidoscopic coloring for $G=K_n$ to be an $(n-3)$-kaleidoscope depicted in the proof of Lemma \ref{lemma2.2} here. When $n$ is even, for each $i\ (1\leq i\leq \frac{n}{2}-2)$, we give a proper coloring to $H_i$ with the colors $2i-1$ and $2i$. Furthermore, we assign the color $n-3$ to all edges in $F$. As for $H_{\frac{n}{2}-1}$ containing vertices $v_1,v_2,v_3,\cdots,v_n$ in the clockwise order. We assign the color $i\ (1\leq i\leq \frac{n}{2})$ to the two edges incident with $v_{2i}$ in $H_{\frac{n}{2}-1}$. While when $n\ (n\geq7)$ is odd, let $v\in V(G)$, then $G-v$ can be decomposed into $\frac{n-1}{2}-1$ Hamiltonian cycles $H_1,H_2,H_3,\cdots,H_{\frac{n-1}{2}-1}$ and a perfect matching $F$, then color the edges of $H_1,H_2,H_3,\cdots,H_{\frac{n-1}{2}-2}$ and $F$ as above. For $H_{\frac{n-1}{2}-1}$ containing $v_1,v_2,\cdots,v_{n-1}$ in the clockwise order, assign the color $i\ (1\leq i\leq \frac{n-1}{2})$ to the edge $v_{2i-1}v_{2i}$ and the color $n-3$ to all the rest edges. At last give the color $i\ (1\leq i\leq \frac{n-1}{2})$ to the edge $vv_{2i-1}$, assign the color $\frac{n-1}{2}+i\ (1\leq i\leq \frac{n-1}{2}-3)$ to the edge $vv_{2i}$ and give the color $n-3$ for the rest of edges in $G$. We denote the coloring depicted above by $c$ and we give a $k$-kaleidoscopic coloring $c'$ on the foundation of $c$ for $\lceil\frac{n}{2}\rceil\leq k\leq n-4$. That is:
\begin{eqnarray*}
c'(e)=\left\{
\begin{array}{rcl}
c(e)   &      & if\ 1\leq c(e)\leq k-1,\\
k      &      & if\ k\leq c(e)\leq n-3.\\
\end{array} \right.
\end{eqnarray*}
This solution can easily be observed since the $[k]$-tuples for distinct vertices are different in their former $\lceil\frac{n}{2}\rceil-1$ positions, so we can still distinguish them after we combine the latter colors into one.

\textbf{Case 2.} For $4\leq k\leq \lceil\frac{n}{2}\rceil-1$, we again distinguish two subcases according to $n$ is even or not.

\textbf{Subcase 1.} When $n$ is even. We first consider the case when $5\leq k\leq \lceil\frac{n}{2}\rceil-1$. Separate the $n$ vertices into two $K_{\frac{n}{2}}$s denoted by $G_1$ and $G_2$ respectively with all edges between them. According to our induction, there exists an $(k-2)$-kaleidoscopic coloring for $K_{\frac{n}{2}}$. Then we give this coloring using the colors from $\{3,4,\cdots,k\}$ to both $G_1$ and $G_2$ with vertices $v_1,v_2,\cdots,v_{\frac{n}{2}}$ for the former and $v'_1,v'_2,\cdots,v'_{\frac{n}{2}}$ for the latter. The vertices $v_i$ and $v'_i\ (1\leq i\leq \frac{n}{2})$ share the same $\{3,4,\cdots,k\}$-tuple. Thus we only need to consider a coloring for the graph $K_{\frac{n}{2},\frac{n}{2}}$ using the colors from $\{1,2\}$ such that each vertex is incident to at least one edge of each color from $\{1,2\}$, besides $v_i$ and $v'_i$ can't have the same $\{1,2\}$-tuple. Let $w_i=v_{i-1}\ (2\leq i\leq \frac{n}{2})$, $w_1=v_{\frac{n}{2}}$, $w'_i=v'_{i+1}\ (1\leq i\leq \frac{n}{2}-1)$ and  $w'_{\frac{n}{2}}=v'_{1}$. The edge $w_iw'_j$ has the color $1$ if and only if $\frac{n}{2}+1\leq i+j\leq n-1$. The color $2$ is assigned to all the remaining edges in $K_{\frac{n}{2},\frac{n}{2}}$. In this way, $v_1$ with $\{1,2\}$-tuple $(2,\frac{n}{2}-2)$ and $v'_1$ with $\{1,2\}$-tuple $(\frac{n}{2}-1,1)$, $v_{\frac{n}{2}-1}$ with $\{1,2\}$-tuple $(\frac{n}{2}-1,1)$ and $v'_{\frac{n}{2}-1}$ with $\{1,2\}$-tuple $(\frac{n}{2}-2,2)$, $v_{\frac{n}{2}}$ with $\{1,2\}$-tuple $(1,\frac{n}{2}-1)$ and $v'_{\frac{n}{2}}$ with $\{1,2\}$-tuple $(\frac{n}{2}-1,1)$, $v_i~(2\leq i\leq \frac{n}{2}-2)$ with $\{1,2\}$-tuple $(i+1,\frac{n}{2}-i-1)$ and $v'_i~(2\leq i\leq \frac{n}{2}-2)$ with $\{1,2\}$-tuple $(i-1,\frac{n}{2}-i+1)$ have different $\{1,2\}$-tuples. As a result, for any $1\leq i\leq \frac{n}{2}$, $v_i$ and $v'_i$ have distinct $\{1,2\}$-tuples and then this coloring is just what we want.

If $k=4$, then $G_i\ (i=1,2)$ contains the unique connected spanning subgraph $F_i$ with only two vertices sharing the same degree as we say in the introduction. Display the vertices of $F_i$ according to their degrees in the nondecreasing order and label them as $v_1,v_2,\cdots,v_{\frac{n}{2}}$ in $F_1$ and $v'_1,v'_2,\cdots,v'_{\frac{n}{2}}$ in $F_2$. Let $H_1=F_1-v_{\frac{n}{2}}v_{\lfloor\frac{n}{4}\rfloor}$ and $H_2=F_2-v'_{\frac{n}{2}}v'_{\lfloor\frac{n}{4}\rfloor}$. Assign the color $3$ to edges of $H_i\ (i=1,2)$ and the color $4$ to the remaining edges in $G_i$s. The rest edges are colored as the above description except for a little change that the color of $w'_{\frac{n}{2}-2}w_{\frac{n}{2}}$ is $2$ instead of $1$. The checkout is similar.

\textbf{Subcase 2.} When $n$ is odd. The simplest condition is $k=4$ and $n\geq13$. In this case, split the graph $K_n$ into $K_{\lceil\frac{n}{2}\rceil}$ and $K_{\lfloor\frac{n}{2}\rfloor}$ with all edges between them. Give a $3$-kaleidoscopic coloring respectively to $K_{\lceil\frac{n}{2}\rceil}$ and $K_{\lfloor\frac{n}{2}\rfloor}$ using the colors $2$, $3$ and $4$. And assign the color $1$ to all edges between them. For the case that $5\leq k\leq \lceil\frac{n}{2}\rceil-1$ and $n\geq13$, label the vertices of $K_{\lceil\frac{n}{2}\rceil}$ and $K_{\lfloor\frac{n}{2}\rfloor}$ respectively as $v_1,v_2,\cdots,v_{\lceil\frac{n}{2}\rceil}$ and $v'_1,v'_2,\cdots,v'_{\lfloor\frac{n}{2}\rfloor}$, furthermore we each give them a $(k-3)$-kaleidoscopic coloring using the colors from $\{4,5,\cdots,k\}$ if $k\neq5$ and each give them a $3$-kaleidoscopic coloring using colors from $\{3,4,5\}$ if $k=5$. For the former case, $v_iv'_j$ is assigned to the color $a\ (a\in\{1,2,3\})$ if $i+j\equiv a\ (\text{mod} \ 3)$ while in the latter case, $v_iv'_j$ has the color $b\ (b\in\{1,2\})$ if $i+j\equiv b\ (\text{mod} \ 2)$.

\begin{figure}[h,t,b,p]
\begin{center}
\includegraphics[scale=0.9]{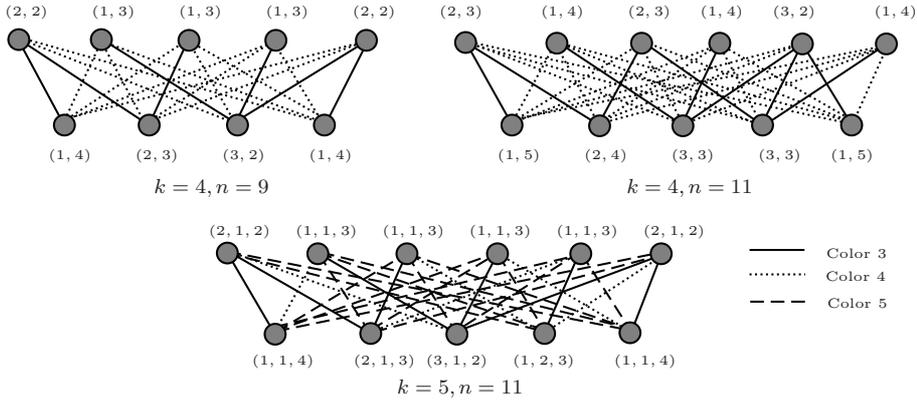}
\caption{Three colorings for $K_{4,5}$ and $K_{5,6}$ in three different conditions}\label{fig1}
\end{center}
\end{figure}

Only three particular conditions are left. That is, $k=4$ when $n=9$ and $k=4$ or $5$ when $n=11$. Like our discussion in Subcase 1, we again find the three unique connected spanning subgraph $F_4$, $F_5$ and $F_6$ with only two vertices of the same degree contained in $K_4$, $K_5$ and $K_6$ appearing in the decomposition of $K_9$ and $K_{11}$. Similarly, the vertices of $F_4$, $F_5$ and $F_6$ are ordered in nondecreasing sequence according to their degrees as $v_i\ (1\leq i\leq4)$, $v'_i\ (1\leq i\leq5)$ and $v''_i\ (1\leq i\leq6)$. Let $H_4=F_4-v_2v_4$, $H_5=F_5-v'_2v'_5$ and $H_6=F_6-v''_3v''_6$. Assign the edges in $H_i\ (i=4,5,6)$ with color the $1$ and all the remaining edges in $K_i\ (i=4,5,6)$ the color $2$. Thus $v_1$ and $v_2$, $v_3$ and $v_4$, $v'_1$ and $v'_2$, $v'_4$ and $v'_5$, $v''_2$ and $v''_3$, in addition with $v''_5$ and $v''_6$, these $6$ couples have the same $\{1,2\}$-tuples. As a result, we only need to provide a coloring to the edges of $K_{4,5}$ or $K_{5,6}$ using the colors $3,4$ and to the edges of $K_{5,6}$ using the colors $3,4,5$ such that the above six couples can not have the same $\{3,4\}$-tuples or $\{3,4,5\}$-tuples. And it goes without saying that the vertices in Figure \ref{fig1} can be matched properly to $v_i\ (1\leq i\leq4)$, $v'_k\ (1\leq k\leq5)$ or $v''_t\ (1\leq t\leq6)$ so that the six couples would not have the same multiset-colors. Thus, our proof is done. \qed

\section{Proof of Conjecture \ref{conjecture1.2}}

\begin{figure}[h,t,b,p]
\begin{center}
\includegraphics[scale=0.9]{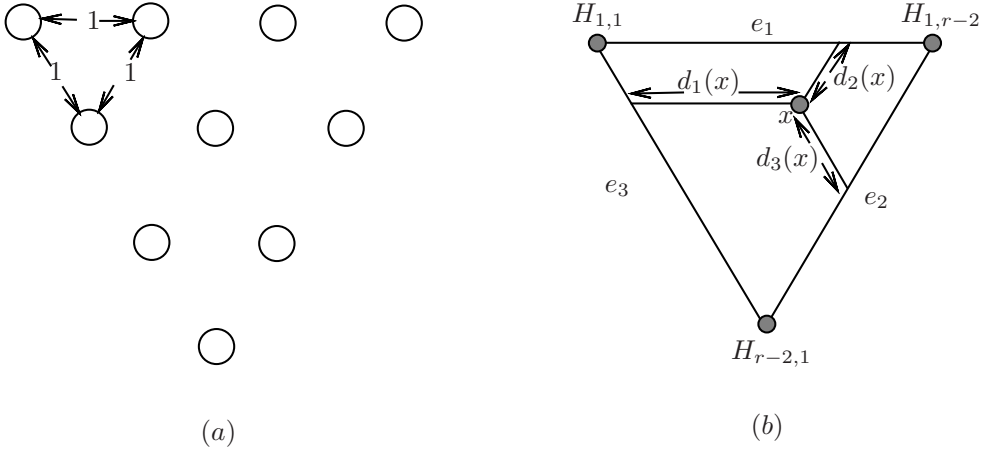}
\caption{The location of the vertices for $r=6$ and the diagram\protect\\ for $d_i(x)$s (i=1,2,3)}\label{fig2}
\end{center}
\end{figure}

For a fixed integer $r\geq5$, note that $\binom{r-1}{2}=\sum_{i=1}^{i=r-2}i$. So for $1\leq i\leq r-2$, let $H_i$ be a set of $r-1-i$ vertices, and for each $H_i$, order the vertices as $H_{i,1},H_{i,2},\cdots,H_{i,r-1-i}$. Then we arrange these vertices in $H_i\ (1\leq i\leq r-2)$ in the shape of an equilateral triangle. That is, the distance between any couple of nearest vertices is $1$. We put an example of the location of the vertices for $r=6$ in Figure \ref{fig2}$(a)$. We then rotate the triangle around the center vertex anticlockwise through an angle of $2\pi/3$ and denote like above as $H'_i\ (1\leq i\leq r-2)$ and $H'_{i,j}\ (1\leq j\leq r-1-i)$. In the same way $H''_i\ (1\leq i\leq r-2)$ and $H''_{i,j}\ (1\leq j\leq r-1-i)$ are obtained after a rotating of $4\pi/3$. We denote by $e_1$ the edge $H_{1,1}H_{1,r-2}$, $e_2$ the edge $H_{1,r-2}H_{r-2,1}$, and $e_3$ the edge $H_{1,1}H_{r-2,1}$. For any point $x$ inside of the triangle, we denote by $d_i(x)\ (i=1,2,3)$ the distance from $x$ to $e_{i+2}\ (\text{mod} \ 3)$ on a line segment parallel to $e_i$. More details are showed below in Figure \ref{fig2}$(b)$. It is obvious that for any vertex $x$, $d_1(x)+d_2(x)+d_3(x)=d(e_1)=r-3$, that is, the length of $e_1$. However, any vertex $x$ can be denoted as $H_{i,s_1(x)}$, $H'_{j,s_2(x)}$ or $H''_{k,s_3(x)}$ where $s_i(x)=d_i(x)+1$. As a result $s_1(x)+s_2(x)+s_3(x)=r$, thus every vertex is endowed with a unique coordinate $(s_1(x),s_2(x),s_3(x))$.
We are going to give a coloring that enables any vertex $x$ to have $s_i(x)\ (i=1,2,3)$ edges with the color $i$.

Firstly we remove the vertex $H_{r-2,1}$. For each $H_i\ (1\leq i\leq r-3)$, we construct the unique connected graph with only two vertices of the same degree satisfying that $\text{deg}_{F_i}H_{i,1}\leq \text{deg}_{F_i}H_{i,2}\leq\cdots\leq \text{deg}_{F_i}H_{i,r-1-i}$. And put $A=\{H_{1,\frac{r-1}{2}+2i}H_{2,\frac{r-1}{2}+2i-1}:1\leq i\leq \frac{r-3}{4}\}$. For $1\leq j\leq \frac{r-3}{4}$, $E_j=\{H_{i,r-i-2j}H_{i+2,r-i-2j-1}:1\leq i\leq r-4j-1\}$. And $B=\{H_{4i-3,\frac{r-4i+3}{2}}H_{4i-1,\frac{r-4i+1}{2}}:1\leq i\leq \frac{r-3}{4}\}$. Note that $A,B$ and $E_j$s are all independent edge sets. The edges in $A,B$, $E_j$s and $F_i$s are all given the color $1$. Then for any vertex $x$, it is incident to exactly $s_1(x)$ edges with the color $1$. The edge sets $F'_i$s in $H'_is\ (2\leq i\leq r-2)$, $F'_1$ in $H'_1-H'_{1,r-2}$, $A'$, $E'_j$s are defined similarly. However $B'=\{H'_{4i-1,\frac{r-4i+1}{2}}H'_{4i+1,\frac{r-4i-1}{2}}:1\leq i\leq \frac{r-3}{4}\}$. Then any vertex $x$ has $s_2(x)$ neighboring edges with the color $2$ except $H'_{2,r-3}$ and $H'_{1,\frac{r-1}{2}}$ with a difference of $1$ respectively. So we add an edge $H'_{2,r-3}H'_{1,\frac{r-1}{2}}$ and give it the color $2$.

Since the vertex $H''_{1,1}$ has been removed, take the vertex $H''_{2,1}$ in place of $H''_{1,1}$ ($H''_{2,1}$ is used twice, both in $F''_{1}$ and $F''_{2}$). Then $F''_i$s in $H''_i$s, $A''$, $E''_j$s and $B''$ are obtained after the same procedure as $F_i$s in $H_i$s, $A$, $E_j$s and $B'$. The number of the edges with the color $3$ incident to $H''_{1,\frac{r-1}{2}}$ is $s_3(H''_{1,\frac{r-1}{2}})-1$ and the number of $H''_{2,1}$ is $s_3(H''_{2,1})+1$. So remove the edge $H''_{2,1}H''_{2,r-3}$ and add the edge $H''_{1,\frac{r-1}{2}}H''_{2,r-3}$ with the color $3$. Then any vertex $x$ is incident to $s_3(x)$ edges with the color $3$. We finally denote by $W_i\ (i=1,2,3)$ the set of edges colored with $i$.

\begin{figure}[h,t,b,p]
\begin{center}
\includegraphics[scale=0.9]{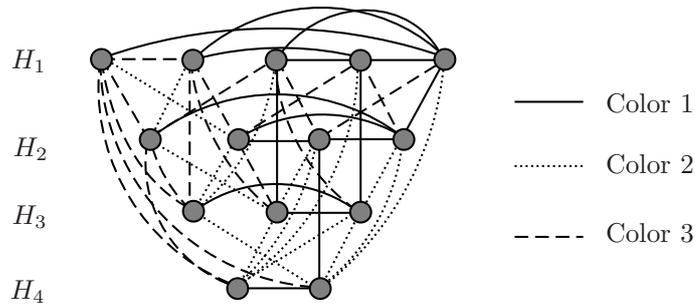}
\caption{A $3$-kaleidoscopic coloring for the $r$-regular graph when $r=7$ }\label{fig3}
\end{center}
\end{figure}

The construction for the $r$-regular graph is complete and the coloring is given. It can be easily verified that this construction satisfies our demands. Besides, any edge in this graph belongs to exactly one of $W_1,W_2,W_3$. We give an example for $r=7$ as Figure \ref{fig3} in the above.

The proof is thus complete.  \qed

\end{document}